\documentclass[12pt]{article}
\usepackage{amssymb,p7}
\newcommand{\PSbox}[3]{\mbox{\rule{0in}{#3}\includegraphics{#1}\hspace{#2}}} 
\newtheorem{Theorem}{Theorem}[section]
\newtheorem{Definition}{Definition}[section]

\newtheorem{Lemma}[Theorem]{Lemma}
\newtheorem{Proposition}[Theorem]{Proposition}

\newtheorem{Corollary}[Theorem]{Corollary}

\newtheorem{Fact}{Fact}[section]
\newtheorem{Remark}{Remark}[section]
\newcommand{\lra}{\longrightarrow}
\def\sqr#1#2{{\vcenter{\vbox{\hrule  height.#2pt
       \hbox{\vrule width.#2pt height#1pt \kern#1pt \vrule width.#2pt}
        \hrule height.#2pt}}}}
\def\bb{\sqr66}   
\let\epsilon=\varepsilon

\let\phi=\varphi
\def\){ \right) }
\def\({ \left( }
\def\[{ \left[ }
\def\]{ \right] }
\def\<{ \langle }
\def\>{ \rangle }
\let\ljunk=\{
\let\rjunk=\}
\def\{{\left\ljunk}
\def\}{\right\rjunk}
\def\p{\partial}

\def\Riem{{\cal R}{\mathrm i}{\mathrm e}{\mathrm m}}

\def\Kod{{\mathrm K}{\mathrm o}{\mathrm d}}

\newcommand{\YW}{{\mathsf Y}\!{\mathsf W}}

\def\lra{\longrightarrow}

\def\Vol{{\mathrm V}{\mathrm o}{\mathrm l}}
\def\Ric{\mbox{Ric}}

\def\ov{\overline}

\newcommand{\R}{{\mathbf R}}

\newcommand{\cpt}{{\mathrm c}{\mathrm p}{\mathrm t}}

\newcommand{\C}{{\mathbf C}}

\def\P{\mbox{\bf P}}
\newenvironment{lproof}[1]{\par\vspace{2mm} \noindent{\bf
#1.} }{\hfill $\bb$\par\vspace{2mm}}
\begin{document} 
\title{The Weyl functional near the Yamabe invariant} 
\author{Kazuo
Akutagawa, Boris Botvinnik, Osamu Kobayashi, Harish Seshadri} 
\date{ \ } 
\maketitle
\markboth{The Weyl functional near the Yamabe invariant}{{\small
\underline{Akutagawa, Botvinnik, Kobayashi, Seshadri, 
The Weyl functional near the Yamabe invariant}}}
\vspace{-14mm}

\begin{abstract}
\noindent
For a compact manifold $M$ of $\dim M =n \geq 4$, we study two conformal
invariants of a conformal class $C$ on $M$. These are the Yamabe
constant $Y_C(M)$ and the $L^{\frac{n}{2}}$-norm $W_C(M)$ of the Weyl
curvature. We prove that for any manifold $M$ there exists a conformal
class $C$ such that the Yamabe constant $Y_C(M)$ is arbitrarily close
to the Yamabe invariant $Y(M)$, and, at the same time, the constant
$W_C(M)$ is arbitrarily large. We study the image of the map $\YW:
C\mapsto (Y_C(M),W_C(M))\in \R^2$ near the line $\{(Y(M),w) \ | \ w\in
\R\}$. We also apply our results to certain classes of
$4$-manifolds, in particular, minimal compact K\"ahler surfaces of Kodaira
dimension $0$, $1$ or $2$.
\end{abstract}
\vspace{-7mm}

\section{Introduction: results and examples}\label{int}
{\bf \ref{int}.1. The Yamabe constant/invariant.} Let $M$ be a smooth
compact (without boundary) manifold of $\dim M =n \geq 3$. We denote by
$\Riem(M)$ the space of the Riemannian metrics on $M$, and by ${\cal
C}(M)$ the space of conformal classes of Riemannian metrics on
$M$. The Einstein-Hilbert functional $I: \Riem(M) \to \R$ is given as
$$
I(g)= \frac{\int_M R_g d\sigma_g}{\Vol_g(M)^{\frac{n-2}{n}}},
$$
where $R_g$ is the scalar curvature and $d\sigma_g$ is the volume
form of $g$. It is well-known that the functional $I$ is not bounded from
above and below for any manifold, and the set of critical points of $I$
coincides with the Einstein metrics on $M$. Let $C\in {\cal C}(M)$ be
a conformal class. The restiction $I|_C$ is always bounded from
below. The constant
$$
Y_C(M):= \inf_{g\in C} I(g)
$$
is known as the \emph{Yamabe constant} of the
conformal class $C$. The Yamabe constant satisfies the inequality
$Y_C(M)\leq Y_{C_0}(S^n)$, where the equality holds if and only if the
manifold $(M,C)$ is conformally equivalent to the standard sphere
$S^n$ with the standard conformal class $C_0$. Then the supremum
$$
Y(M):= \sup_{C\in {\cal C}(M)} Y_C(M)
$$
is the \emph{Yamabe invariant} of $M$.  The Yamabe constant $Y_C(M)$
is an  important conformal invariant. In particular, for any
conformal class $C\in {\mathcal C}(M)$, there exists a metric (Yamabe
metric) $\check{g}\in C$ such that $I(\check{g})=Y_C(M)$. The Yamabe
metric has constant scalar curvature $R_{\check{g}}=Y_C(M)$ under the
normalization $\Vol_{\check{g}}(M)=1$.
\vspace{2mm}

\noindent
{\bf \ref{int}.2. The Weyl functional.} Now let $\dim M\geq 4$. We
denote by $W_g=(W^{i}_{\ jk\ell})$ the Weyl tensor of a metric $g\in
\Riem(M)$. The norm $|W_g|_g$ is defined as
$$
|W_g|_g = \left(W^{i}_{\ jk\ell} W_{i}^{\ jk\ell}\right)^{1/2}.
$$
Let $C\in {\cal C}(M)$, and let $g\in C$ be any metric. Then the
integral 
$$
W_C(M):= \int_M |W_g|^{\frac{n}{2}}_g d\sigma_g
$$
does not depend on the choice of the metric $g\in C$. The map ${\cal
C}(M)\to \R$, $C\mapsto W_C(M)$, is called the \emph{Weyl functional}. We
shall call $W_C(M)$ the \emph{Weyl constant} of the conformal class
$C$ and the infimum
$$
W(M):=\inf_{C\in {\cal C}(M)} W_C(M)
$$
the \emph{Weyl invariant} of $M$. Clearly the invariant $W(M)$ is a
diffeomorphism invariant as well as $Y(M)$.
\vspace{2mm}

\noindent
{\bf \ref{int}.3. General results.} To state our main results, we
define one more numerical invariant $\omega(M)$ as follows. Firstly,
for a given manifold $M$, we say that a sequence $\{ C_i\}$, $C_i\in
{\mathcal C}(M)$, of conformal classes is a \emph{Yamabe sequence} if
$$
\lim_{i\to\infty} Y_{C_i}(M) =  Y(M).
$$
In particular, if there exists a class $C\in {\mathcal C}(M)$
satisfying $Y(M)=Y_C(M)$, the sequence $C_i=C$ for $i\geq 1$ is a
Yamabe sequence. Secondly, we define the constant
$$
\omega(M,\{ C_i\}):=\displaystyle \!\!\!\!\!\!  \lim_{\ \ \ \
i\to\infty}\!\!\!\!\!\inf W_{C_i}(M)\geq 0
$$
for any Yamabe sequence $\{ C_i\}$. Then the invariant $\omega(M)$ is
defined as
$$
\omega(M) :=\displaystyle
\inf \{ \omega(M,\{ C_i\}) \ | \ \mbox{$\{ C_i\}$ is a Yamabe sequence}\}.
$$
Notice that $\omega(M)\geq W(M)$ by definition. Also it easy to show
that $\omega(M)=0$ for $M= S^n$, $S^1\times S^{n-1}$, $T^n$ and their
connected sums.
\vspace{2mm}

\noindent
{\bf Theorem A.} {\sl Let $M$ be a compact manifold of $\dim M = n\geq
4$. Then for any small $\epsilon>0$ and any constant $\kappa
>\omega(M)$ there exists a conformal class $C\in {\mathcal C}(M)$ such
that}
\begin{equation}\label{int_1}
\{
\begin{array}{l}
Y_C(M)\geq Y(M)-\epsilon,
\\
\kappa +\epsilon \geq W_C(M) \geq  \kappa .
\end{array}
\right.
\end{equation}
\begin{Remark}
{\rm If $\omega(M)=+\infty$, Theorem A delivers an empty statement. It
is not clear that the invariant $\omega(M)$ is finite for every compact
manifold $M$. However, without the finiteness  of
$\omega(M)$, the following statement still holds.}
\end{Remark}
{\bf Theorem A$^{\prime}$.} {\sl Let $M$ be a compact manifold
of $\dim M = n\geq 4$. Then for any small $\epsilon>0$ and any
constant $\kappa >0$ there exists a conformal class $C\in {\mathcal
C}(M)$ such that}
$$
\{
\begin{array}{l}
Y_C(M)\geq Y(M)-\epsilon,
\\
W_C(M) \geq  \kappa .
\end{array}
\right.
$$
We prove Theorem A in two steps. Firstly, we prove Theorem A for
$M=S^n$ by constructing a conformal class $C\in {\mathcal C}(S^n)$
satisfying (\ref{int_1}) with the understanding that
$\omega(S^n)=0$. Secondly, we prove the general case by constructing
an appropriate conformal class on the connected sum $M\#
S^n$. Theorem A$^{\prime}$ follows immediately from our argument as
well.
\vspace{2mm}

\noindent
To proceed further, we would like to introduce some terminology.
\vspace{2mm}

\noindent
{\bf \ref{int}.4. The $\YW$-quadrant, Yamabe corner, Sobolev and
Kuiper lines.}  For a compact manifold $M$ we consider the map
$$
\YW : {\cal C}(M) \lra \R^2, \ \ \ C\mapsto (Y_C(M), W_C(M)).
$$
We denote by ${\mathbf K}^{\YW}(M)\subset \R^2$ the image of the map
$\YW$. By definition, the set ${\mathbf K}^{\YW}(M)$ is a 
diffeomorphism invariant of $M$.
\vspace{2mm}

\noindent
The third author observed that the set ${\mathbf K}^{\YW}(M)$ contains
certain interesting aspects of conformal geometry, and he studied some
of its properties (cf. \cite{Kobayashi1a}, \cite{Kobayashi0}).
\vspace{2mm}

\noindent
Let $(y,w)$ be the Euclidian coordinates in $\R^2$. Consider the
 {\emph{$\YW$-quadrant} of a manifold $M$:
$$
{\mathbf Q}^{\YW}(M):=\{ (y,w) \ | \ y\leq Y(M), \ \ w\geq W(M) \}, 
$$
see Fig. \ref{int}.1. Clearly ${\mathbf K}^{\YW}(M)\subset {\mathbf
Q}^{\YW}(M)$. We emphasize that {\bf there is not much known about the
shape of the set ${\mathbf K}^{\YW}(M)\subset {\mathbf
Q}^{\YW}(M)$.} In this paper we study the set ${\mathbf
K}^{\YW}(M)$ near the Sobolev line defined below.
\vspace{2mm}

\noindent
\hspace*{54mm}\PSbox{sob1.pstex}{5cm}{4cm}
\begin{picture}(0,0)
\put(-167,72){{\small (the Kuiper line)}}
\put(-127,32){{\small $w=W(M)$}}
\put(-75,110){{\small $w$}}
\put(-50,-3){{\small $Y(M)$}}
\put(3,62){{\small $y=Y(S^n)$}}
\put(3,43){{\small $\omega(M)\geq W(M)$}}
\put(18,25){{\small (the Yamabe corner)}}
\put(0,80){{\small $y=Y(M)$ (the Sobolev line)}}
\put(15,5){{\small $y$}}
\end{picture}
\vspace{3mm}

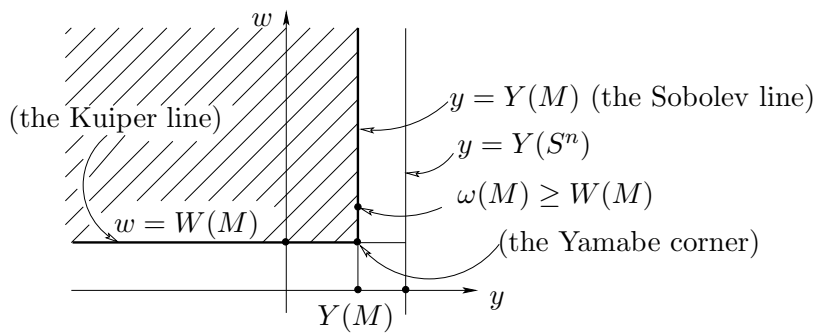
\centerline{{\small {\bf Fig. \ref{int}.1.} The $\YW$-quadrant of $M$.
Here $Y(M)>0$, $W(M)>0$.}}
\vspace{2mm}

\noindent
We explain the notations and terminology given in Fig. \ref{int}.1. We
call the point $(Y(M),W(M))$ the \emph{Yamabe corner} of $M$.  Then we
call the line $\{y=Y(M)\}$ the \emph{Sobolev line}.  This name is
motivated by the following observation. Let $Y(M)>0$, and let $C$ be a
positive conformal class (i.e. $Y_C(M)>0$). We denote by $\check{g}\in
C$ a Yamabe metric. Then the following inequality holds
(cf. \cite{Akutagawa1})
$$
\left( \int_M |f|^{\frac{2n}{n-2}} d\sigma_{\check{g}}
\right)^{\frac{n-2}{n}} \leq \frac{1}{\frac{n-2}{4(n-1)}Y_C(M)}\int_M
|df|^2_{\check{g}} d\sigma_{\check{g}}
+\frac{1}{\Vol_{\check{g}}(M)^{\frac{2}{n}}}\int_M f^2
d\sigma_{\check{g}}
$$
for $f\in L^{1,2}_{\check{g}}(M)$, where $L^{1,2}_{\check{g}}(M)$
denotes the Sobolev space of $L^2$-functions with $L^2$ first
derivatives relative to $\check{g}$. In particular, the constant
$\frac{n-2}{4(n-1)}Y_C(M)$ is the best Sobolev constant for the
Sobolev embedding $ L^{1,2}_{\check{g}}(M) \subset
L^{\frac{2n}{n-2}}_{\check{g}}(M) $. Furthermore,
$\frac{n-2}{4(n-1)}Y(M)$ is the supremum of those Sobolev constants.
We notice that
$$
{\mathbf K}^{\YW}(S^n)\cap \{\mbox{the Sobolev line of $S^n$}\}
= \{(Y(S^n),0)\} =\{\mbox{the Yamabe corner of $S^n$}\}
$$
for all $n\geq 4$ due to the final resolution of the Yamabe problem by
Schoen \cite{Schoen0}. The case $Y(M)\leq 0$ does not have such
interpretation, however, it is convenient to use this term in general
case.
\vspace{2mm}

\noindent
We call the line $\{w=W(M)\}$ the \emph{Kuiper line} of $M$ (see
Gursky's paper \cite{Gursky1} for some remarkable results on the Kuiper
lines of $4$-manifolds).  Our motivation for this name comes from the
Kuiper Theorem \cite{Kuiper}: if $W_C(S^n)=0$ for a conformal class
$C\in {\cal C}(S^n)$, then $C$ is equivalent (up to a diffeomorphism)
to the standard conformal class $C_0$. In our terms,
$$
{\mathbf K}^{\YW}(S^n)\cap  \{\mbox{the Kuiper line of $S^n$}\} 
= \{(Y(S^n),0)\} =\{\mbox{the Yamabe corner of $S^n$}\}.
$$
In terms of the set ${\mathbf K}^{\YW}(M)$, Theorem A
implies the following.
\vspace{2mm}

\noindent
{\bf Corollary B.} {\sl Let $M$ be a compact manifold of $\dim M
= n\geq 4$. Then
$$
\ov{{\mathbf K}^{\YW}(M)} \cap \{\mbox{the Sobolev line of $M$}\}=
\{y = Y(M), \ w\geq \omega(M)\}.
$$
Here $\ov{{\mathbf K}^{\YW}(M)}$ is the closure of ${\mathbf
K}^{\YW}(M)\subset \R^2$.}
\vspace{1mm}

\noindent
{\bf \ref{int}.5. The Einstein and the gap curves.} Here we consider
only 4-dimensional manifolds. For a manifold $M$ we denote by
$\chi(M)$ and $\tau(M)$ the Euler characteristic and the signature of
$M$ respectively.  The Hirzebruch signature formula gives the
following result.
\begin{Proposition}\label{Hirz} 
Let $M$ be an oriented, compact $4$-manifold. Then for any
$C\in {\mathcal C}(M)$, 
$$
W_C(M)\geq 48\pi^2|\tau(M)|
$$
with the equality if and only if $(M,C)$ is half conformally flat
(i.e. self-dual or anti self-dual). In particular, $ W(M)\geq
48\pi^2|\tau(M)|$.
\end{Proposition}
The Chern-Gauss-Bonnet formula leads to the
following interesting result.
\begin{Theorem}\label{Gauss-Bonnet} {\rm (\cite[Section 1-(i)]{Kobayashi0})}
Let $M$ be a compact $4$-manifold. Then for any conformal class $C\in
{\mathcal C}(M)$ the following inequality holds:
$$
W_C(M) \geq 32 \pi^2 \chi(M) -\frac{1}{6}Y_C(M)^2 
$$
with the equality if and only if the conformal class $C$ contains an
Einstein metric.
\end{Theorem}
We call the curve $ w = 32 \pi^2 \chi(M) -\frac{1}{6}y^2$, $w\geq
W(M)$, $y\leq Y(M)$, the \emph{Einstein curve} of $M$.
\vspace{2mm}

\noindent
On the other hand, the Gap Theorems due to Gursky \cite[Theorem
3.3, Proposition 3.5]{Gursky4} and Gursky--LeBrun \cite[Theorem
1]{Gursky3} can be reformulated as follows.
\begin{Theorem}\label{gap} {\rm (Gap Theorem)}
Let $M$ be an oriented, compact $4$-manifold and $C\in {\mathcal
C}(M)$ a conformal class satisfying $Y_C(M)>0$.
Assume that either
\begin{enumerate}
\item[{\bf (i)}] $b^+_2(M) = \dim H_+^2(M;\R)>0$, or
\item[{\bf (ii)}] $C$ contains an Einstein metric  and $W_g^+\not\equiv 0$
for $g\in C$.
\end{enumerate}
Then the following inequality holds:
$$
W_C(M) \geq \frac{1}{3}Y_C(M)^2  - 48\pi^2\tau(M).
$$
Here $W_g^+$ is the self-dual Weyl curvature of $g$.
\end{Theorem}
We call the curve $w=\frac{1}{3} y^2- 48\pi^2\tau(M)$, $w\geq 0$,
$y\geq 0$, the \emph{gap curve} of $M$.
\vspace{2mm}

\noindent
{\bf \ref{int}.6. Examples. 1.} Let $M=S^4$. Then we have $W(S^4)=0$,
$Y(S^4)= 8\pi\sqrt{6}$. The Einstein and the gap curves are $ w=
64\pi^2 - \frac{1}{6} y^2$, $ w =\frac{1}{3} y^2 $ respectively.  We
notice that $Y(S^4)=8\pi\sqrt{6}$ coincides with the value of $y$
given by the intersection of the Einstein curve with $w=0$.  We have
$$
\begin{array}{lcl}
{\mathbf K}^{\YW}(S^4)\cap\{\mbox{the Sobolev line of $S^4$}\} &=&
\{\mbox{the Yamabe corner of $S^4$}\}, 
\\ 
\\ 
{\mathbf
K}^{\YW}(S^4)\cap\{\mbox{the Kuiper line of $S^4$}\} &=& \{\mbox{the
Yamabe corner of $S^4$}\}.
\end{array}
$$
On the other hand, Corollary B implies that
$$
\ov{{\mathbf K}^{\YW}(S^4)}\cap \{\mbox{the Sobolev line of $S^4$}\}
= \{y = 8\pi\sqrt{6}, \ w \geq 0\}. 
$$
We have the following picture for $S^4$:
\vspace{2mm}

\noindent
\hspace*{54mm}\PSbox{sob2.pstex}{5cm}{53mm}
\begin{picture}(0,0)
\put(-142,87){{\small (the Einstein curve)}}
\put(-30,145){{\small $w$}}
\put(45,-7){{\small $Y(S^4)=8\pi\sqrt{6}$}}
\put(10,2){{\small $8\pi\sqrt{2}$}}
\put(85,15){{\small $y$}}
\put(-40,52){{\small $64\pi^2$}}
\put(-42,34){{\small $\frac{128}{3}\pi^2$}}
\put(55,35){{\small (the gap curve)}}
\put(55,120){{\small $128\pi^2$}}
\end{picture}
\vspace{3mm}

\centerline{{\small {\bf Fig. \ref{int}.2.} The $\YW$-picture of $S^4$.}}
\vspace{2mm}

\noindent
Here and below we shade a subset of ${\mathbf Q}^{\YW}(M)$ which still
contains ${\mathbf K}^{\YW}(M)$. We call it the ``$\YW$-picture of
$M$''.
\vspace{2mm}

\noindent
Then Theorem \ref{gap} implies that the intersection $ {\mathbf
K}^{\YW}(S^4)\cap\{\mbox{the Einstein curve of $S^4$}\} $ does not
contain the points with $ 8\pi\sqrt{2}< y < 8\pi\sqrt{6}.  $ It is not
known whether this intersection contains any points except the Yamabe
corner.
\vspace{2mm}

\noindent
{\bf Example 2.} Let $M=\C\P^2$. The Einstein curve and the gap curves
are given as $ w = 98 \pi^2- \frac{1}{6}y^2$ and $w= \frac{1}{3}y^2- 48\pi^2$
respectively.  
\vspace{2mm}

\noindent
It is known that $W(\C\P^2)=48\pi^2$ (it follows from Proposition
\ref{Hirz}) and $Y(\C\P^2)=12\pi\sqrt{2}$ (LeBrun \cite{LeBrun6}, 
Gursky--LeBrun \cite{Gursky2}), where $Y(\C\P^2)$ is attained by the
conformal class $[g_{FS}]$ of the Fubini-Study metric $g_{FS}$ on
$\C\P^2$. In particular,
$\omega(\C\P^2)=W_{[g_{FS}]}(\C\P^2)=W(\C\P^2)= 48\pi^2$.  The
$\YW$-picture of $\C\P^2$ is shown at Fig. \ref{int}.3 below.  Here
the Yamabe corner coincides with the intersection of the Sobolev line,
Kuiper line and the Einstein curve (and the gap curve). Moreover, it
follows from Gursky--LeBrun's result
\cite[Theorem 7]{Gursky2} that
$$
\begin{array}{c}
\!\!{\mathbf K}^{\YW}(\C\P^2)\!\cap\!
\{\mbox{the Sobolev line of $\C\P^2$}\}\!=\!
\{\mbox{the Yamabe corner of $\C\P^2$}\}\! =\!\{(12\pi\sqrt{2}, \!48\pi^2)\}.
\end{array}
$$
Corollary B gives here that
$$
\begin{array}{c}
\!\!
\ov{{\mathbf K}^{\YW}(\C\P^2)}\!\cap\! 
\{\mbox{the Sobolev line of $\C\P^2$}\}\! =\!
\{ y= 12\pi\sqrt{2}, \ w \geq 48\pi^2\}.
\end{array}
$$
In this case the gap curve does not give useful restrictions.
However, it follows from \cite[Theorem 2]{Gursky3} that if a conformal
class $C$ contains an Einstein metric different from $g_{FS}$, then
$Y_C(\C\P^2)< 4 \pi\sqrt{6}$, see Fig \ref{int}.3.
\vspace{2mm}

\noindent
\hspace*{54mm}\PSbox{sob3a.pstex}{5cm}{53mm}
\begin{picture}(0,0)
\put(-152,107){{\small (the Einstein curve)}}
\put(-55,154){{\small $w$}}
\put(-70,74){{\small $96\pi^2$}}
\put(-70,32){{\small $48\pi^2$}}
\put(-43,22){{\small $4\pi\!\sqrt{6}$}}
\put(-45,-7){{\small $Y(\C\P^2)\!=\!12\pi\sqrt{2}$}}
\put(40,30){{\small (the gap curve)}}
\put(35,10){{\small $y$}}
\end{picture}
\vspace{3mm}

\centerline{{\small {\bf Fig. \ref{int}.3.} The $\YW$-picture of $\C\P^2$.}}

\vspace{2mm}

\noindent
{\bf \ref{int}.7. Results on 4-manifolds.} 
Let $\chi(M)$ and
$\tau(M)$ be as above. 
\vspace{2mm}

\noindent
{\bf Theorem C.} {\sl Let $M$ be a minimal compact complex surface
of general type, that is, its Kodaira dimension $\Kod(M)=2$. Let
$M^{\prime}= M\# k\ov{\C\P^2} \#\ell(S^1\times S^3)$ ($k,\ell\geq 0$)
be the connected sum of the blow-up of $M$ at $k$ points with $\ell$
copies of $S^1\times S^3$. Then} 
$$
\begin{array}{c}
\omega(M^{\prime}) =\frac{16}{3}\pi^2(4\chi(M^{\prime})- 
3 \tau(M^{\prime})+2k+8\ell), \ \  \ 
\mbox{in particular,}
\\
\\
\ov{{\mathbf K}^{\YW}(M^{\prime})}\!\cap \!
\{\mbox{{\sl the Sobolev line of} $M^{\prime}$}\}\! =\!
\{ y\!= \!Y(M^{\prime}), \ w \geq \frac{16}{3}\pi^2(4\chi(M^{\prime})\!- \!
3 \tau(M^{\prime})+2k+8\ell)\}.
\end{array}
$$
In this case the Yamabe invariant $Y(M^{\prime})$ is known (see LeBrun
\cite{LeBrun5}, \cite{LeBrun2} and Petean \cite{Petean}):
$$
Y(M^{\prime}) = -4\sqrt{2}\pi\sqrt{2\chi(M) + 3\tau(M)} < 0.
$$
We have the following $\YW$-picture for such $M$:
\vspace{2mm}

\noindent
\hspace*{54mm}\PSbox{sob4b.pstex}{5cm}{55mm}
\begin{picture}(0,0)
\put(5,88){{\small $\left(\!\!\begin{array}{c}\mbox{the Einstein curve}\\
w = 32\pi^2\chi(M) - \frac{1}{6}y^2\end{array}\!\!\right)$}}
\put(-37,140){{\small $w$}}
\put(-15,32){{\small $\omega(M)=\frac{16}{3}\pi^2(4\chi(M)- 3 \tau(M))$}}
\put(-16,18){{\small $48\pi^2|\tau(M)|$}}
\put(-125,-10){{\small $Y(M)\!=\!-4\sqrt{2}\pi\sqrt{2\chi(M) + 3\tau(M)}$}}
\put(70,10){{\small $y$}}
\put(-135,25){{\small ? ? ? ?}}
\end{picture}
\vspace{3mm}

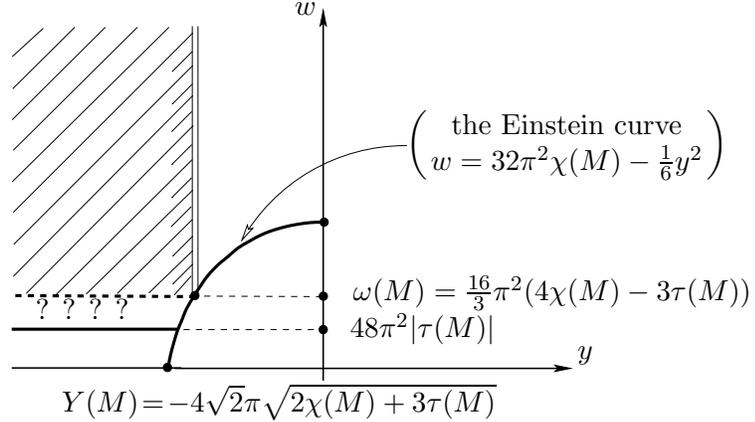
\centerline{{\small {\bf Fig. \ref{int}.4.} The $\YW$-picture of a
minimal compact complex surface $M$ of general type.}}
\vspace{2mm}

\noindent
In Fig. \ref{int}.4, we assume that $M$ has a negative first Chern
class.  We emphasize that in general, it is not clear whether
$\omega(M)=W(M)$ or $\omega(M)>W(M)$ for every minimal compact complex
surface of general type.
\vspace{2mm}

\noindent
\noindent
{\bf Theorem D.} {\sl Let $M$ be a minimal compact K\"ahler-type
complex surface of $\Kod(M)=0$ or $\Kod(M)=1$. Let $M^{\prime}= M\#
k\ov{\C\P^2} \#\ell(S^1\times S^3)$ ($k,\ell\geq 0$). Then}
$$
\begin{array}{c}
\omega(M^{\prime}) = W(M^{\prime})= -48\pi^2\tau(M^{\prime}) = 
48\pi^2(k-\tau(M)), \ \  \ 
\mbox{in particular,}
\\
\\
\ov{{\mathbf K}^{\YW}(M^{\prime})}\cap 
\{\mbox{{\sl the Sobolev line of} $M^{\prime}$}\} =
\{ y= Y(M^{\prime}), \ w \geq -48\pi^2\tau(M^{\prime})\}.
\end{array}
$$
In this case it is also known that $Y(M^{\prime})=0$ (see LeBrun
\cite{LeBrun2} and Petean \cite{Petean}).
\begin{Remark}{\rm {\bf (1)}
Under the assumptions of Theorem C or D, the intersection
$$
{\mathbf K}^{\YW}(M) \cap \{\mbox{the Sobolev line of $M$}\} = \{\mbox{one
point}\} \ \ \mbox{or} \ \ \emptyset ,
$$
see \cite[Proposition 5.89]{Aubin} or \cite[Theorem 1]{Kobayashi}.  }
\vspace{1mm}

\noindent
{\bf (2)} {\rm Theorems C and D still hold for more general $4$-manifolds
$M^{\prime}$, see \cite[Theorems A, B and C]{Ishida}.}
\end{Remark}
\bf \ref{int}.8. More examples.} Here we give examples illustrating
Theorems C and D.
\vspace{2mm}

\noindent
{\bf Example 3.} Let $M$ be a $K3$-surface. Then $M$ is a minimal
K\"ahler surface of $\Kod(M)=0$.  Here $\tau(M)=-16$, and Theorem D
gives that $W(M)=\omega(M)=768\pi^2$, $Y(M)=0$.  Here we know also
that there is no conformal class $C\in {\mathcal C}(M)$ with
$Y_C(M)=0$ and with $W_C(M)> \omega(M)=W(M)= 768\pi^2$. Thus we have:
$$
\begin{array}{c}
{\mathbf K}^{\YW}(M)\cap \{\mbox{the Sobolev line}\} = \{(Y(M)=0, 
W(M)= 768\pi^2)\} = \{\mbox{the Yamabe corner}\}, 
\\
\\
\ov{{\mathbf K}^{\YW}(M)}\cap \{\mbox{the Sobolev line}\} =
\{ y= 0 , \ w \geq  768\pi^2\} \ \ \ \mbox{by Corollary B.}
\end{array}
$$
In particular, we see that a conformal class $C\in {\mathcal C}(M)$
contains an Einstein metric if and only if $Y_C(M)=0$ and
$W_C(M)=768\pi^2$. We have the following $\YW$-picture:

\noindent
\hspace*{54mm}\PSbox{sob5b.pstex}{5cm}{52mm}
\begin{picture}(0,0)
\put(7,32){{\small $\left(\!\!\begin{array}{c}\mbox{the Einstein curve}\\
w = 768\pi^2 - \frac{1}{6}y^2\end{array}\!\!\right)$}}
\put(-37,140){{\small $w$}}
\put(-15,60){{\small $\omega(M)=W(M)= 768\pi^2$}}
\put(-40,-10){{\small $Y(M)=0$}}
\put(70,10){{\small $y$}}
\end{picture}
\vspace{3mm}

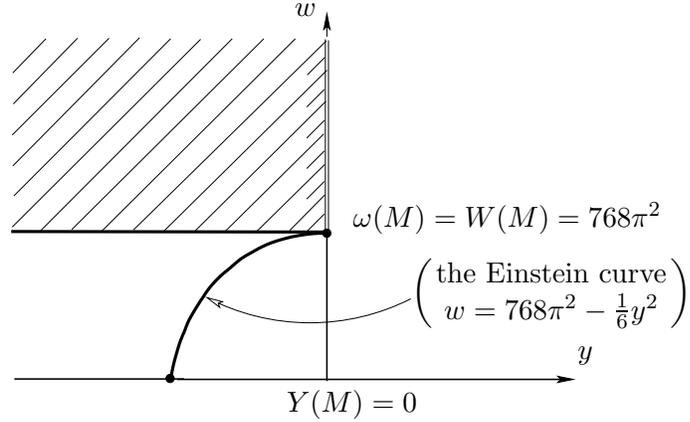
\centerline{{\small {\bf Fig. \ref{int}.5.} The $\YW$-picture of the
$K3$ surface.}}
\vspace{1mm}

\noindent
{\bf Example 4.} Let $M=T^2\times \Sigma_g$, where $\Sigma_g$ is a
surface with genus $g\geq 2$. Then $M$ is a minimal K\"ahler surface of
$\Kod(M)=1$. Here $\tau(M)=0$ and $\chi(M)=0$.
\vspace{2mm}

\noindent
Theorem D gives that $\omega(M)=W(M)=0$. We notice that $Y(M)=0$,
however, there is no conformal class $C\in {\mathcal C}(M)$ such that
$Y(M)=Y_C(M)=0$ (see \cite[Proposition 6]{LeBrun2}).  Hence we have
$$
\begin{array}{c}
{\mathbf K}^{\YW}(M)\cap \{\mbox{the Sobolev line of $M$}\} = \emptyset, 
\\
\\
\ov{{\mathbf K}^{\YW}(M)}\cap \{\mbox{the Sobolev line of $M$}\} =
\{ y=0, \ w \geq  0 \} \ \ \ \mbox{by Corollary B.}
\end{array}
$$
\noindent
\hspace*{54mm}\PSbox{sob6b.pstex}{50mm}{60mm}
\begin{picture}(0,0)
\put(65,27){{\small $\left(\!\!\begin{array}{c}\mbox{the Einstein curve}\\
w = - \frac{1}{6}y^2\end{array}\!\!\right)$}}
\put(-60,150){{\small $w$}}
\put(-45,35){{\small $\omega(M)=W(M)= 0$}}
\put(-70,10){{\small $Y(M)=0$}}
\put(30,25){{\small $y$}}
\end{picture}
\vspace{3mm}

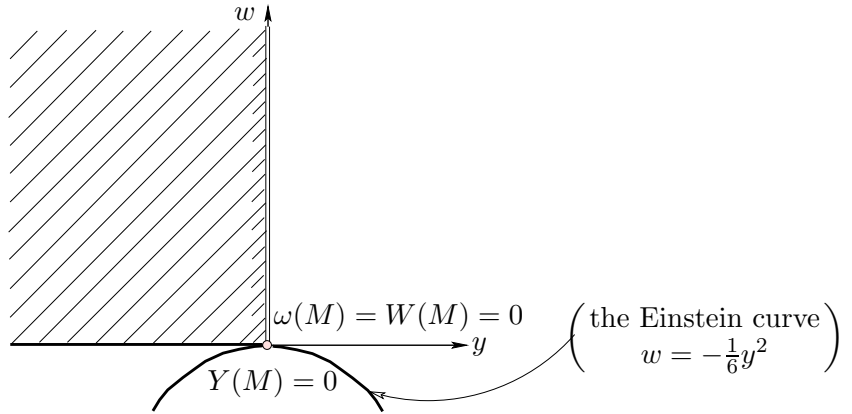
\centerline{{\small {\bf Fig. \ref{int}.6.} The $\YW$-picture of 
$M=T^2\times \Sigma_g$, $g\geq 2$.}}
\vspace{2mm}

\noindent
{\bf Example 5.} Let $M=\Sigma_{g_1}\times\Sigma_{g_2}$ where $g_1, \
g_2\geq 2$. Then $M$ is a minimal K\"ahler surface of $\Kod(M)=2$. We
have $\chi(M)= 4(g_1-1)(g_2-1)$, $\tau(M)=0$.  Here the Yamabe
invariant is attained by the canonical product K\"ahler-Einstein
metrics.  Theorem C gives that
$\omega(M)=\frac{256}{3}\pi^2(g_1-1)(g_2-1)$. Also $Y(M)=
-16\pi\sqrt{(g_1-1)(g_2-1)}$. We have
$$
\begin{array}{c}
{\mathbf K}^{\YW}(M)\cap \{\mbox{the Sobolev line of $M$}\} = 
\{(Y(M), \omega(M))\} = \{\mbox{the Yamabe corner of $M$}\},
\\
\\
\ov{{\mathbf K}^{\YW}(M)}\cap \{\mbox{the Sobolev line of $M$}\} =
\{ y\geq Y(M), \  w \geq  \omega(M) \} \ \ \ \mbox{by Corollary B.}
\end{array}
$$
Here it is not clear whether $\omega(\Sigma_{g_1}\times\Sigma_{g_2})>
W(\Sigma_{g_1}\times\Sigma_{g_2})$ or they are equal.
\vspace{2mm}

\noindent
\hspace*{40mm}\PSbox{sob7b.pstex}{5cm}{55mm}
\begin{picture}(0,0)
\put(5,88){{\small $\left(\!\!\begin{array}{c}\mbox{the Einstein curve}\\
w = 128\pi^2 (g_1-1)(g_2-1)- \frac{1}{6}y^2\end{array}\!\!\right)$}}
\put(-37,140){{\small $w$}}
\put(-15,30){{\small $\omega(M)=\frac{256}{3}\pi^2(g_1-1)(g_2-1)$}}
\put(-125,-10){{\small $Y(M)\!=\!-16\pi\sqrt{(g_1-1)(g_2-1)}$}}
\put(70,10){{\small $y$}}
\put(-135,15){{\small ? ? ? ?}}
\end{picture}
\vspace{4mm}

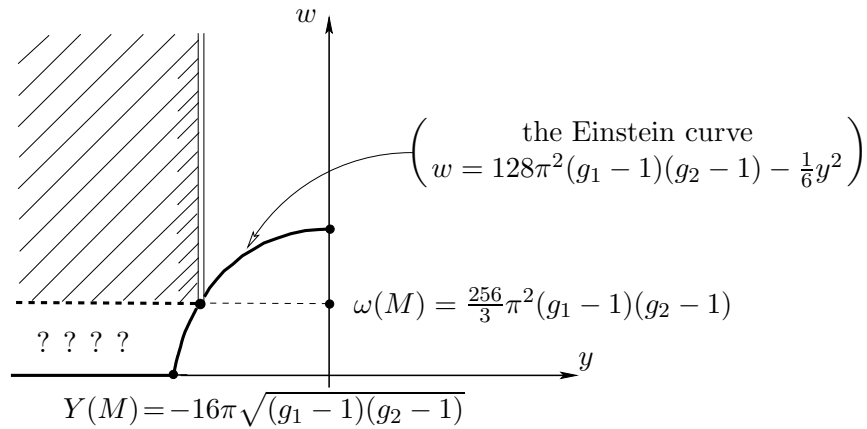
\centerline{{\small {\bf Fig. \ref{int}.7.} The $\YW$-picture of
$M=\Sigma_{g_1}\times\Sigma_{g_2}$ with $g_1, \ g_2\geq 2$.}}
\vspace{2mm}

\noindent
{\bf Example 6.} Let $M=\C H^2/\Gamma$ be a smooth compact quotient of
the complex hyperbolic space. Then $M$ is a minimal compact K\"ahler
surface of general type. Then $\chi(M)=3\tau(M)> 0$ by \cite[Theorem
5]{LeBrun3}.
\vspace{2mm}

\noindent
Theorem D now gives:
$
\omega(M)=\frac{16}{3}\pi^2\cdot 9\tau(M)= 48\pi^2\tau(M) = 16\pi^2\chi(M).
$ Let $g_B$ be the canonical K\"ahler-Einstein metric on $M$ (i.e. the
Bergmann metric). The Yamabe invariant $Y(M)=
-12\sqrt{2}\pi\sqrt{\tau(M)}$ is attained by the metric $g_B$. Since
$g_B$ is a self-dual metric, we have
$W(M)= 48\pi^2\tau(M) $ by Proposition \ref{Hirz}. Hence $\omega(M)=W(M)=
48\pi^2\tau(M) = 16\pi^2\chi(M)$.  Here we have
$$
\begin{array}{c}
{\mathbf K}^{\YW}(M)\cap \{\mbox{the Sobolev line of $M$}\} = 
\{(Y(M), \omega(M))\} = \{\mbox{the Yamabe corner of $M$}\},
\\
\\
\ov{{\mathbf K}^{\YW}(M)}\cap \{\mbox{the Sobolev line of $M$}\} =
\{ y\geq Y(M), \  w \geq  \omega(M)\} 
\ \ \ \mbox{by Corollary B.}
\end{array}
$$

\noindent
\hspace*{50mm}\PSbox{sob8b.pstex}{5cm}{50mm}
\begin{picture}(0,0)
\put(5,88){{\small $\left(\!\!\begin{array}{c}\mbox{the Einstein curve}\\
w = 96 \tau(M)- \frac{1}{6}y^2\end{array}\!\!\right)$}}
\put(-37,140){{\small $w$}}
\put(-15,30){{\small $\omega(M)=W(M)=48\tau(M)$}}
\put(-125,-10){{\small $Y(M)\!=\!-12\sqrt{2}\pi\sqrt{\tau(M)}$}}
\put(70,10){{\small $y$}}
\end{picture}
\vspace{4mm}

\centerline{{\small {\bf Fig. \ref{int}.8.} The $\YW$-picture of
$M=\C H^2/\Gamma$.}}
\vspace{2mm}

\noindent
{\bf \ref{int}.9. The plan of the paper.} Firstly we prove Theorems C,
D assuming Theorem A in Section \ref{s1}. Secondly we prove Theorem A
for the particular manifold $M=S^n$ in Section \ref{s2}. Thirdly we prove
the general case of Theorem A in Section \ref{s3}.
\vspace{2mm}

\noindent
{\bf \ref{int}.10. Acknowledgments.}  The first and the third authors
are very grateful to the Department of Mathematics of the University
of Oregon for hospitality. The authors are grateful to Claude LeBrun
for his valuable comments concerning Theorem \ref{gap} and Theorems C,
D.

\section{Proofs of Theorems C and D}\label{s1}
{\bf \ref{s1}.1. Proof of Theorem C.} Let $M$ be a minimal compact
complex surface of general type, and $M^{\prime}=M\# k\ov{\C\P^2}
\#\ell(S^1\times S^3)$. Here we use the results by LeBrun
\cite[Theorem 7]{LeBrun4}, \cite[Theorem 2]{LeBrun2},  
\cite[Proposition 3.2, Proof
of Theorem 4.3]{LeBrun1} and Petean \cite[Proposition 4]{Petean} to
conclude that there exists a sequence of metrics $\{g_j\}$, $g_j\in
\Riem(M^{\prime})$ such that
$$
\{
\begin{array}{l}
\displaystyle 
\lim_{j\to\infty}\int_{M^{\prime}} |W_{g_j}^+|^2_{g_j}d\sigma_{g_j}
= \frac{4}{3} \pi^2(2\chi(M)+3\tau(M))>0, 
\\ 
\\ 
\displaystyle 
\lim_{j\to\infty} Y_{[g_j]}(M^{\prime}) = Y(M^{\prime})
= -4\sqrt{2}\pi\sqrt{2\chi(M)+3\tau(M)}<0,
\\ 
\\ 
\displaystyle 
\lim_{j\to\infty} \int_{M^{\prime}} R^2_{g_j} d\sigma_{g_j}=  
Y(M^{\prime})^2 = 
32\pi^2(2\chi(M)+3\tau(M)) >0.
\end{array}
\right.
$$
Here we use the following convention: $|W_g|_g^2= 4\left(|W_g^+|_g^2 +
|W_g^-|_g^2
\right)$. 
 In particular, the sequence of conformal classes $\{ [g_j] \}$ is a
Yamabe sequence. Then the Hirzebruch signature formula gives:
$$
\begin{array}{rcl}
\displaystyle \lim_{j\to\infty} W_{[g_j]}(M^{\prime}) 
&=& \displaystyle  
8\lim_{j\to\infty}\int_{M^{\prime}} 
|W_{g_j}^+|^2_{g_j}d\sigma_{g_j} -48\pi^2\tau(M^{\prime})
\\
\\
&=&
\frac{16}{3}\pi^2(4\chi(M^{\prime})-3\tau(M^{\prime})+ 2k+8\ell).
\end{array}
$$
Then, by definition, we have: 
$$
\omega(M^{\prime},\{[g_j]\}) = \lim_{j\to\infty} W_{[g_j]}(M^{\prime})
= \frac{16}{3}\pi^2(4\chi(M^{\prime})-3\tau(M^{\prime})+ 2k+8\ell).
$$
From the formula (\ref{LB}) below, we also notice
that
$$
\omega(M^{\prime}) \geq \frac{16}{3}\pi^2(4\chi(M^{\prime})-3\tau(M^{\prime})+ 2k+8\ell).
$$
These imply Theorem C.  \hfill $\Box$
\begin{Remark}
{\rm Under the same assumptions on $M$ and $M^{\prime}$ 
as in Theorem C, we notice the following.
\begin{enumerate}
\item[{\bf (1)}] There is a sharp inequality
\begin{equation}\label{LB}
|Y_C(M^{\prime})| +\sqrt{6}\left(
\int_{M^{\prime}} |W_g^+|_g^2 d\sigma_g
\right)^{\frac{1}{2}}\geq 6\sqrt{2}\pi \sqrt{2\chi(M)+3\chi\tau(M)}>0
\end{equation}
for any $C\in {\mathcal C}(M^{\prime})$ and $g\in C$ (see
\cite[Formula (13)]{LeBrun1}). When either $k\geq 1$ or $\ell\geq 1$,
the inequality (\ref{LB}) provides a further restriction on the
$\YW$-picture of $M^{\prime}$.
\item[{\bf (2)}] Furthemore,
$$
\begin{array}{c}
\{\mbox{the Sobolev line of $M^{\prime}$}\}\cap\{\mbox{the Einstein curve of
$M^{\prime}$}\} \ \ \ \ \ \ \ \ \ \ \ \ \ \ \ \ \ \ \ \ \ \ \ \ 
\\
\\
\ \ \ \ \ \ \ \ \ \ \ \ 
= \{(Y(M^{\prime}),\frac{16}{3}\pi^2(4\chi(M^{\prime}) -
3\tau(M^{\prime}) -k -4\ell))\}.
\end{array}
$$
Then if either $k\geq 1$ or $\ell\geq 1$, the above $w$-coordinate is
strictly less than $\omega(M^{\prime})$.
\end{enumerate}}
\end{Remark}
\vspace{1mm}

\noindent
{\bf \ref{s1}.2. Proof of Theorem D.} Let $M$ be a minimal compact
K\"ahler-type complex surface of Kodaira dimension $\Kod(M)=0$ or
$1$. Let $M^{\prime}=M\#k\ov{\C\P^2}\#\ell(S^1\times S^3)$.  Again we
use the results due to LeBrun \cite[Theorem 7]{LeBrun4},
\cite[Theorems 4 and 6]{LeBrun2}, \cite[Proposition 3.2 and Proof of
Theorem 4.3]{LeBrun1} and Petean \cite[Proposition 3]{Petean} to
conclude that there exists a sequence of metrics $\{g_j\}$, $g_j\in
\Riem(M^{\prime})$ such that
$$
\{
\begin{array}{l}
\displaystyle 
\lim_{j\to\infty}\int_{M^{\prime}} |W_{g_j}^+|^2_{g_j}
d\sigma_{g_j}= 
\frac{4}{3}\pi^2(2\chi(M)+3\tau(M)),
\\ 
\\ 
\displaystyle 
\lim_{j\to\infty} Y_{[g_j]}(M^{\prime}) = Y(M^{\prime})
= 0,
\\ 
\\ 
\displaystyle 
\lim_{j\to\infty} \int_{M^{\prime}} R^2_{g_j} d\sigma_{g_j}=  
Y(M^{\prime})^2 = 0.
\end{array}
\right.
$$
Thus the sequence of conformal classes $\{ [g_j] \}$ is a Yamabe
sequence. Then the Hirzebruch signature formula gives:
$$
\begin{array}{rcl}\displaystyle
\lim_{j\to\infty} W_{[g_j]}(M^{\prime}) &=&\displaystyle
8\lim_{j\to\infty}\int_{M^{\prime}} |W_{g_j}^+|^2_{g_j}
d\sigma_{g_j} -48\pi^2\tau(M^{\prime})
\\
\\
 &=&\displaystyle \frac{16}{3}\pi^2(4\chi(M)-3\tau(M)+9k)
\end{array}
$$
We obtain:
$$
\omega(M^{\prime},\{[g_j]\})= \lim_{j\to\infty} 
W_{[g_j]}(M^{\prime}) = \frac{16}{3}\pi^2(4\chi(M)-3\tau(M)+9k).
$$
It also follows from \cite[Theorem 4 and Proposition 7]{LeBrun2} that
$ 2\chi(M)=-3\tau(M)\geq 0 $ and $\omega(M^{\prime})=W(M^{\prime})=
48\pi^2(k-\tau(M))= -48\pi^2\tau(M^{\prime})$.  These imply Theorem
D.\hfill $\Box$

\section{The sphere}\label{s2}
{\bf \ref{s2}.1. The result.} Our goal here is to prove Theorem A for
$M=S^n$, that is, the following result.
\begin{Theorem}\label{s2_Th1}
For any small $ \epsilon > 0$ and arbitrary $\kappa > 0$ there exists a
conformal class $C\in {\mathcal C}(S^n)$ such that
$$
\{
\begin{array}{l}
Y_C (S^n) \ge Y(S^n) - \epsilon, 
\\
\kappa + \epsilon \ge W_C(S^n) \ge \kappa  .
\end{array}
\right.
$$
\end{Theorem}
To give a proof of Theorem \ref{s2_Th1} we will construct a required
conformal class $C\in {\mathcal C}(S^n)$ using the manifold
$S^{n-1}\times \R$. The following definition is due to Schoen and Yau
\cite{SY1}.
\begin{Definition}\label{s2_def1}
{\rm Let $(N,C)$ be a conformal manifold of $\dim N =n \geq 3$ (possibly
noncompact) and $g\in C$ any metric. We denote
$$ 
\begin{array}{rcl}
E_g(f) &=&\displaystyle
 \int _N \left( \alpha_n { \vert df \vert }^2_g + R_g f^2 \right) 
d\sigma_g ,
\\
\\
Q_g (f)&=&\displaystyle
{ E_g (f) \over ( \int _N {\vert f \vert} ^{2n \over n-2} 
d\sigma_g )^{n-2 \over n}}, 
\end{array}
$$
where $\alpha_n= \frac{4(n-1)}{n-2}$. Then the Yamabe constant $Y_{C}
(N)$ is defined as
$$
Y_{C} (N)=\displaystyle \inf_{\begin{array}{c}_{f \in
C_{\cpt}^{\infty}(N)} \\ _{f \not \equiv 0}\end{array}} Q_g(f).
$$ 
The Yamabe constant $Y_{C} (N)$ does not depend on the choice of the
metric $g\in C$ (see \cite{SY1}).}
\end{Definition}
{\bf \ref{s2}.2.  Continuity of the Yamabe constant.} First we recall
the following result on the continuity of the Yamabe constant, see
\cite{Besse}, \cite[Fact 1.4]{Kobayashi} (cf. \cite[Lemma 4.1]{AB1}).
\begin{Fact}\label{s2_prop3}
Let $M$ be a compact manifold  of $\dim M=n\geq 3$ and $g\in
\Riem(M)$. Let $\{ g_i\} $ be a sequence of Riemannian metrics such that
$$
\{
\begin{array}{l}
g_i \to g,
\\
R_{g_i} \to R_{g}
\end{array}
\right.
\ \ \ \ \ \ \mbox{as \ \ \ \ \ $i\to\infty$} 
$$
in the uniform $C^0$-norm on $M$ with respect to $g$. Then
$Y_{[g_i]}(M) \to Y_{[g]}(M)$ as $i\to\infty$. 
\end{Fact}
In the course of proving Proposition \ref{s2_prop3} one has to show
that $\displaystyle \lim_{\ \ \ \ i\to\infty}\!\!\!\!\!\inf
Y_{[g_i]}(M)\geq Y_{[g]}(M)$.  Here it  is essential that the volume
$\Vol_g(M)$ is finite. However, a similar continuity result still
holds for noncompact manifolds under the positivity of scalar
curvature.
\begin{Proposition}\label{s2_prop4}
Let $N$ be a manifold without boundary (possibly noncompact). Assume
that a metric $g\in \Riem(N)$ satisfies the following
condition
\begin{equation}\label{s2_eq2a}
0< L_0^{-1}\leq R_g \leq L_0 \ \ \ \mbox{on} \ \ N
\end{equation}
for some constant $L_0>0$. Let $\{g_i\}$ be a sequence of Riemannian
metrics, such that
$$
\{
\begin{array}{l}
g_i \to g,
\\
R_{g_i} \to R_{g}
\end{array}
\right.
\ \ \ \ \ \ \mbox{as \ \ \ \ $i\to\infty$} 
$$
in the uniform $C^0$-norm on $N$ with respect to $g$. Then
$Y_{[g_i]}(N) \to Y_{[g]}(N)$ as $i\to\infty$.
\end{Proposition}
\begin{Proof}
Given $ \epsilon > 0$, there exists an integer $i(\epsilon)>0$ such
that
$$ 
\begin{array}{c}
\ \ \ \ \ \ \ \ \ \ \ \ \ \ \ \ \ \ \ \ \ \ \ \ \ 
\vert g_i - g \vert_g \le \epsilon, \ \ \ \ \ 
\vert g_i^{-1} - g^{-1} \vert_g \le \epsilon \ \ \ \mbox{(with respect to $g$),}
\\
\\
\begin{array}{rcccl}
 (1- \epsilon )d\sigma_g &\le& d\sigma_{g_i} &\le& (1+ \epsilon )d\sigma_g,
\\
\\
 0 < (1- \epsilon)R_g &\le& R_{g_i} &\le& (1+ \epsilon)R_g
\end{array}
\end{array}
$$
for any $i\geq i( \epsilon)$. Note that the last inequality follows
from the condition (\ref{s2_eq2a}). For any compactly supported smooth
function $ f \in C_{\cpt}^\infty (N) $, we have
$$
\begin{array}{rcl}
 E_{g_i}(f) &=& \displaystyle 
\int _N ( \alpha_n { \vert df \vert}_{g_i} ^2 + R_{g_i} f^2 ) d\sigma_{g_i} 
\\
\\
& \ge & \displaystyle \int _N (\alpha_n(1- \epsilon ) { \vert df \vert}_{g} ^2 +(1-  \epsilon) R_{g} f^2 )(1- \epsilon) d\sigma_{g} 
\\
\\
& \ge & \displaystyle ( 1-2K \epsilon ) 
\int _N ( \alpha_n { \vert df \vert}_{g} ^2 + R_{g} f^2 ) d\sigma_{g}  
\\
\\
& \ge & \displaystyle ( 1-2K \epsilon ) E_g (f)        
\end{array}
$$
for any $i\geq i( \epsilon)$. Similarly we have
$$ 
\begin{array}{c}
E_{g_i} (f) \le (1+2K \epsilon ) 
E_g (f) \ \ \ \mbox{and} 
\\
\\
\displaystyle 
(1- K^{\prime} \epsilon) 
\left( 
\int _N {\vert f \vert} ^{2n \over n-2} d\sigma_{g}   
\right)^{n-2 \over n} 
\le \left( \int _N {\vert f \vert} ^{\frac{2n}{n-2}}
d\sigma_{g_i} 
\right)^{n-2 \over n} \le (1+ K^{\prime}\epsilon) 
\left( \int _N
{\vert f \vert} ^{2n \over n-2} d\sigma_{g} 
\right)^{n-2 \over n} .
\end{array}
$$ 
Here $K$, $K^{\prime}$ are positive constants independent of $f$ and
$\epsilon$.  Hence we have
$$ 
\begin{array}{rcccll}
\displaystyle
{(1-2K \epsilon ) \over (1+ K^{\prime} \epsilon)} Q_g (f) &\le&
Q_{g_i} (f) &\le&\displaystyle
{{(1+2K \epsilon )} \over {(1- K^{\prime}\epsilon)}}
Q_g (f), &\mbox{or}
\\
\\ 
\displaystyle
(1-  K^{\prime\prime} \epsilon) Q_g (f) &\le&  Q_{g_i} (f) &\le&\displaystyle
 (1+ K^{\prime\prime}  \epsilon) Q_g(f) , &\mbox{or}
\\
\\ 
\displaystyle
 (1-  K^{\prime\prime} \epsilon) Y_{[g]} (N) &\le& Y_{[g_i]} (N)
&\le &\displaystyle (1+  K^{\prime\prime} \epsilon) Y_{[g]} (N) &
\end{array}
$$
for some positive constant $K^{\prime\prime}$. This implies that
$Y_{[g_i]} (N)
\longrightarrow Y_{[g]} (N)$.
\end{Proof}
We recall the following well-known facts (cf. \cite{Kobayashi},
\cite{Schoen1}).
\begin{Fact}\label{s2_fact1}
Let $g_{0}$ and $h_{0}$ denote the standard metrics on the spheres
$S^n$ and $S^{n-1}$ respectively. Let $C_0=[g_0]$. Then
$$
Y(S^n)=Y_{C_0}(S^n) = Y_{C_0}(S^n\setminus \{\mbox{2 points}\})=
Y_{[h_0+dt^2]}(S^{n-1}\times \R).
$$
\end{Fact}
\begin{Fact}\label{s2_fact2}
Let $M$ be a compact manifold, $p_1,\ldots,p_k\in M$, and $C\in
{\mathcal C}(M)$ a conformal class which is conformaly flat near
the points $p_1,\ldots,p_k$.  Then
$$
Y_C(M) = Y_C(M\setminus \{p_1,\ldots,p_k\}).
$$
\end{Fact}
Let $g_{0}$ and $h_{0}$, as above, denote the standard metrics on the
spheres $S^n$ and $S^{n-1}$ respectively. As a corollary of
Proposition \ref{s2_prop4} we have the following result.
\begin{Corollary}\label{s2_cor1}
Let $h (\cdot ,t)$ be a smooth (with respect to $t$) family of
Riemannian metrics on $S^{n-1}$, and let $\bar{g} (z,t) = h (z,t) +
dt^2$, $(z,t)\in S^{n-1}\times \R$, be the corresponding metric on the
cylinder $S^{n-1}\times \R$. Then for any $\epsilon>0$ there exists
$\delta>0$ such that if the family $h (\cdot ,t)$ satisfies the
condition
$$
\{
\begin{array}{cl}
|h(\cdot ,t) - h_0|_{h_0} < \delta & 
\mbox{on $S^{n-1}$ for any $t\in \R$, and}
\\
|R_{\bar{g}(\cdot,t)}- R_{\bar{g}_0}| <
\delta & \mbox{on} \ \ S^{n-1}\times \R,
\end{array}
\right.
$$
then $Y_{[\bar{g}]}(S^{n-1}\times \R) \geq
Y_{[\bar{g}_0]}(S^{n-1}\times \R) -\epsilon = Y(S^n) -\epsilon$. Here
$\bar{g}_0:=h_0+dt^2$ on  $S^{n-1}\times \R$.
\end{Corollary}
{\bf \ref{s2}.3. Metrics on the cylinder $N=M\times \R$.} We consider
the following general situation. Let $M$ be a compact manifold, and
$N=M\times \R$.
\vspace{2mm}

\noindent
First, let $g(\cdot,t)\in \Riem(M)$ be a smooth family of metrics and
$\bar{g}(z,t)= g(z,t) +dt^2$ the corresponding metric on $N=M\times
\R$. First, we compare the scalar curvature functions $R_{g}$ and
$R_{\bar{g}}$, where $g:=g(\cdot,t)$. We have:
\begin{equation}\label{s2_eq3a}
\begin{array}{rcl}
R_{\bar{g}} & = & 
R_{g} +\frac{1}{4}
\left(
g^{ij}g^{kl}\cdot \p_t g_{ij}\cdot \p_t g_{kl} - (g^{ij}\cdot \p_t g_{ij})^2
\right)
-\frac{1}{2}
\left(
g^{ij}\cdot \p_t^2g_{ij} + \p_t(g^{ij}\cdot \p_t g_{ij})
\right)
\\
\\
 & = & 
R_{g} -g^{ij}\cdot \p_t^2g_{ij} -\frac{1}{4}(g^{ij}\cdot \p_t g_{ij})^2
+\frac{3}{4} g^{ij}g^{kl}\cdot \p_t g_{ij}\cdot \p_t g_{kl}.
\end{array}
\end{equation}
Here the indices $i,j,k,\ell$ vary within the corresponding indices
$1,\ldots,n-1$ of local coordinates $z=\{z^1,\ldots, z^{n-1}\}$ on
$M$. Now we consider a particular family of metrics. Let $L>0$ be a
sufficiently large constant.  We choose a nonnegative function
$\phi=\phi_L\in C^{\infty}(\R)$ satisfying
$$
|\phi^{\prime}|\leq \frac{2}{L}, \ \ |\phi^{\prime\prime}|\leq \frac{4}{L^2},
\ \ \ \phi= \{ \begin{array}{ll}
0, & \mbox{if} \  t\leq 0,
\\
1, & \mbox{if} \  t\geq L,
\end{array}\ \ \ \mbox{(see Fig. \ref{s2}.1.)}
\right.
$$
\noindent
\hspace*{40mm}\PSbox{sob5.pstex}{5cm}{25mm}
\begin{picture}(0,0)
\put(-65,50){{\small $1$}}
\put(93,0){{\small $L$}}
\end{picture}
\vspace{3mm}

\centerline{{\small {\bf Fig. \ref{s2}.1.} The function $\phi=\phi_L$.}}
\vspace{2mm}

\noindent
Let $h,\hat{h}\in \Riem(M)$ be two metrics on the ``slice'' $M$. We
define the family of metrics on $M\times\{t\}\cong M$:
$$
\begin{array}{rcl}
g(z,t)&:=&\phi(t)\cdot h(z) + (1-\phi(t))\cdot \hat{h}(z)
\\
&=& h(z) -(1-\phi(t))\cdot T(z), 
\end{array}
$$
where $T(z):=h(z) - \hat{h}(z)$. Let
$\bar{g}(z,t):= g(z,t)+ dt^2$ be the metric on the cylinder
$M\times\R$.  We notice that
$$
\{
\begin{array}{rcl}
g^{\prime}&=&\phi^{\prime}\cdot T,
\\
g^{\prime\prime}&=&\phi^{\prime\prime}\cdot T
\end{array}
\right. \ \ \ \mbox{and} \ \ \ \bar{g} = \{
\begin{array}{ll}
\hat{h}+ dt^2 &\mbox{if} \ \ t\leq 0,
\\
h+ dt^2 &\mbox{if} \ \ t\geq L.
\end{array}
\right.
$$
Assume that $|T|_h<<1$, so that $\frac{1}{2}(h_{ij}) \leq (g_{ij}) \leq 
2(h_{ij})$.
Then we use (\ref{s2_eq3a}) to give the estimate
\begin{equation}\label{s2_eq4a}
|R_{\bar{g}}- R_g| \leq \frac{8}{L^2}|T|_h +\frac{4}{L^2}|T|_h^2 +
\frac{12}{L^2}|T|_h^2 \leq \frac{K_0}{L^2}|T|_h 
\end{equation}
since $|T|_h<<1$. We denote $\theta(z,t):= -(1-\phi(t))T(z)$, so that
$g(z,t)= h(z)+ \theta(z,t)$ on $M\times\{t\}\cong M$.  We also denote
$$
P_h(\theta):= -\nabla^j\nabla_i\theta_i^i + \nabla^i\nabla^j \theta_{ij} -
R_{ij}\theta^{ij},
$$
where $(R_{ij})$ is the Ricci curvature of $h$. A straightforward
calculation gives the following formula for the scalar curvature
(cf. \cite{Kobayashi}, \cite{Kobayashi2}):
$$
R_g = R_h + P_h(\theta) + Q_h(\theta), 
$$
where the function $Q_h(\theta)$ is estimated by
\begin{equation}\label{s2_eq5a}
\begin{array}{c}
\displaystyle
|Q_h(\theta)|_h \leq K_n \{ |\nabla\theta|^2_h\cdot q^3 + |\theta|_h 
\cdot|\nabla^2\theta|_h\cdot q^2 + \left(|\theta|_h \cdot |\nabla^2\theta|_h
+|\Ric(h)|_h\cdot |\theta|^2_h\right)q\} \ \ \ \mbox{with}
\\
\\
\displaystyle
q(z,t):=\max\{\left. \frac{h(X,X)}{g(X,X)} \right| 
\ X\in T_{(z,t)} M, \ X\not\equiv 0 \}.
\end{array} \!\!\!\!\!\!\!\!\!\!\!\!\!\!\!
\end{equation}
Here $\nabla$, $\Ric=(R_{ij})$ and all norms are with respect to the
metric $h$. Now we also assume that $|\nabla T|_h$ and $|\nabla^2
T|_h$ are sufficiently small. Then the estimates (\ref{s2_eq4a}) and
(\ref{s2_eq5a}) imply that
$$
|R_g - R_h| \leq \Phi(h,T),
$$
where the constant $\Phi(h,T)\geq 0 $ (depending on $h$ and $T$) is
small as well. We obtain
$$
|R_{\bar{g}} - R_h | \leq \frac{K_0}{L^2}|T|_h + \Phi(h,T).
$$
In particular, the above argument implies the following technical
result.
\begin{Lemma}\label{s2_L7}
Let $h_0\in \Riem(S^{n-1})$ be the standard metric of constant
curvature $1$. Then for any integer $j>0$ there exist a constant
$L(j)>>1$ and a metric $h_j\in \Riem(S^{n-1})$ such that
$$
\{
\begin{array}{l}
|\bar{g}_0 - \bar{g}_j|_{\bar{g}_0} \leq \frac{1}{j},
\\
|R_{\bar{g}_0} - R_{\bar{g}_j}| \leq \frac{1}{j}
\end{array}
\right. \ \ \ \mbox{on} \ \ \ S^{n-1}\times \R.
$$
Here the metrics $\bar{g}_0$ and $\bar{g}_j$ on $S^{n-1}\times \R$ are
defined as
$$
\begin{array}{rcl}
\bar{g}_0(z,t)&:=& h_0 (z) + dt^2,
\\
\bar{g}_j(z,t)&:=& \left(
h_0(z) - (1-\phi(t))\cdot(h_0(z) - h_j(z))
\right) + dt^2,
\end{array}
$$
where $\phi(t):=\phi_L(t)$ for any $L\geq L(j)$.
\end{Lemma}
\begin{lproof}{Proof of Theorem \ref{s2_Th1}}
Let $h_0$ be the standard metric on $S^{n-1}$. 
We use Lemma
\ref{s2_L7} to conclude the following.
\vspace{2mm}

\noindent
For any $\epsilon>0$ there exist a constant $L>>1$ and a metric $h\in
\Riem(S^{n-1})$ (where $h$ is not homothetic but $C^2$-close to $h_0$)
such that the metric
$$
\begin{array}{rcl}
\bar{g}(z,t) & :=& g(z,t) + dt^2 \ \ \mbox{on $S^{n-1}\times\R$ with}
\\
\\
g(z,t) & :=& \!\!\! 
\{
\begin{array}{ll}
h(z) & \mbox{on} \ S^{n-1}\times [-\bar{L}, \bar{L}],
\\
h_0(z) -(1-\phi_L(t-\bar{L}))\cdot (h_0(z) -h(z)) 
 & \mbox{on} \ S^{n-1}\times [\bar{L}, \bar{L}+L],
\\
h_0(z) -(1-\phi_L(-t-\bar{L}))\cdot (h_0(z) -h(z)) 
 & \mbox{on} \ S^{n-1}\times [-(\bar{L}+L), -\bar{L}],
\\
h_0(z) & \mbox{on} \ S^{n-1}\times \left(\R\setminus 
[-(\bar{L}+L), \bar{L}+L]\right)
\end{array}
\right.
\end{array}
$$
satisfies the inequalities
\begin{equation}\label{s2_eq6a}
\{
\begin{array}{rcl}\displaystyle
Y_{[\bar{g}]}(S^{n-1}\times \R)&\geq &Y(S^n)-\epsilon,
\\
\\
\displaystyle
\int_{S^{n-1}\times \left( [-(\bar{L}+L), -\bar{L}] 
\cup [\bar{L}, \bar{L}+L]\right)} 
|W_{\bar{g}}|^{\frac{n}{2}}_{\bar{g}} d\sigma_{\bar{g}}&\leq & \epsilon
\end{array}
\right.
\end{equation}
for any $\bar{L}>0$ (see Fig. \ref{s2}.2).

\noindent
\hspace*{30mm}\PSbox{sob6.pstex}{5cm}{26mm}
\begin{picture}(0,0)
\put(-60,62){{\small $L$}}
\put(80,62){{\small $L$}}
\put(-42,-10){{\small $-\bar{L}$}}
\put(52,-10){{\small $\bar{L}$}}
\put(-132,22){{\small $h_0+dt^2$}}
\put(-10,22){{\small $h+dt^2$}}
\put(120,22){{\small $h_0+dt^2$}}
\end{picture}
\vspace{5mm}

\centerline{{\small {\bf Fig. \ref{s2}.2.} The metric $\bar{g}$ on
$S^{n-1}\times\R$.}}
\vspace{2mm}

\noindent
Notice that the restriction of the metric $\bar{g}(z,t)$ on
$S^{n-1}\times [-\bar{L},\bar{L}]$ is given as $\bar{g}(z,t)= h(z) +
dt^2$.  Consider the function
$$
f(\bar{L}):= \int_{S^{n-1}\times [-\bar{L},\bar{L}]}
|W_{\bar{g}}|^{\frac{n}{2}}_{\bar{g}} d\sigma_{\bar{g}} \in (0,\infty).
$$
Then $f$ is continuous on $(0,\infty)$, and
\begin{equation}\label{s2_eq7a}
\lim_{\bar{L}\to 0} f(\bar{L}) = 0, \ \ \ \ 
\displaystyle
\lim_{\bar{L}\to \infty} f(\bar{L}) = \infty.
\end{equation}
By definition of $f$, we have:
\begin{equation}\label{s2_eq8a}
\int_{S^{n-1}\times \R} |W_{\bar{g}}|^{\frac{n}{2}}_{\bar{g}} 
d\sigma_{\bar{g}} 
= f(\bar{L}) + 
\int_{S^{n-1}\times\left([-(\bar{L}+L),-\bar{L}]
\cup [\bar{L},\bar{L}+L]\right)}
|W_{\bar{g}}|^{\frac{n}{2}}_{\bar{g}}  d\sigma_{\bar{g}} .
\end{equation}
It follows from (\ref{s2_eq6a}), (\ref{s2_eq7a}) and (\ref{s2_eq8a})
that for any $0<\epsilon<< 1$ and any constant $\kappa>0$ there exist
constants $L>0$ and $\bar{L}>0$ such that
\begin{equation}\label{s2_eq9a}
\{
\begin{array}{l}
\bar{g}(z,t)= g(z,t)+ dt^2 \ \ \ \mbox{on $S^{n-1}\times \R$ (as above)},
\\
Y_{[\bar{g}]}(S^{n-1}\times \R) \geq Y(S^n) -\epsilon,
\\
\kappa+\epsilon \geq W_{[\bar{g}]}(S^{n-1}\times \R) \geq \kappa,
\end{array}
\right.
\end{equation}
where $\bar{g}(z,t)= h_0(z) + dt^2$ on $S^{n-1}\times
\left(\R\setminus [-(\bar{L}+L),\bar{L}+L]\right)$. Thus the conformal
class $[\bar{g}]$ can be extended smoothly to a conformal class $C\in
{\mathcal C}(S^n)$ such that
\begin{equation}\label{s2_eq10a}
\{
\begin{array}{rcl}
Y_{[\bar{g}]}(S^{n-1}\times \R) &=& Y_C(S^{n}),
\\
W_{[\bar{g}]}(S^{n-1}\times \R) &=& W_C(S^n).
\end{array}
\right.
\end{equation}
Combining (\ref{s2_eq9a}) with (\ref{s2_eq10a}), we complete the proof
of Theorem \ref{s2_Th1}.
\end{lproof}

\noindent
\hspace*{55mm}\PSbox{sob7.pstex}{5cm}{32mm}
\begin{picture}(0,0)
\put(-10,40){{\small $p_+$}}
\put(-160,40){{\small $p_-$}}
\end{picture}
\vspace{5mm}

\centerline{{\small {\bf Fig. \ref{s2}.3.} The sphere $(S^n, \hat{g})$.}}
\vspace{2mm}

\begin{Remark}
{\rm In the proof of Theorem \ref{s2_Th1}, we let $\bar{L}$ go to the
infinity. Then (in the terminology of \cite{AB2}) we obtain the
canonical cylindrical manifold
$$
\left(
S^{n-1}\times \R, \bar{g}:=h+dt^2\right), 
$$
where $h$ is not homothetic but $C^2$-close to $h_0$.  Proposition
\ref{s2_prop4}  implies that the inequality $
Y_{[\bar{g}]}(S^{n-1}\times \R) \geq Y(S^n)-\epsilon $ still holds for
a small $\epsilon>0$. However, $W_{[\bar{g}]}(S^{n-1}\times
\R)=+\infty$.  From
\cite[Theorems 6.1, 6.2]{AB2} there exists a function $u\in
C^{\infty}_+(S^{n-1}\times
\R)\cap L^{1,2}_{\bar{g}}(S^{n-1}\times \R)$ such that the metric
$\hat{g}= u^{\frac{4}{n-2}} \bar{g}$ is a Yamabe metric (i.e. a
minimizer of the functional $Q_{\bar{g}}$, see \cite{AB2} for
details).  Furthermore, the metric $\hat{g}$ is \emph{almost conical
metric} near the two ends.
\vspace{2mm}

\noindent
Then, by adding two points $p_-$, $p_+$ to
the ends, we obtain the sphere
$$
S^n= (S^{n-1}\times \R)\cup \{p_-,p_+\}
$$
with the metric $\hat{g}$ (see Fig. \ref{s2}.3).
\vspace{2mm}

\noindent
Thus $(S^n, \hat{g})$
is a compact Riemannian manifold with two almost conical singularities.
These singularities provide the source of the phenomenon that
$$
W_{[\hat{g}]}(S^n) = W_{[\bar{g}]}(S^{n-1}\times \R) = + \infty.
$$ 
}
\end{Remark}

\section{Proof of Theorem A}\label{s3}
We start with the following approximation result in 
\cite{Kobayashi} and \cite{Kobayashi2} (see also \cite{AB1}). 
\begin{Proposition}\label{s3_Prop1}
Let $M$ be a compact manifold of $\dim M = n\geq 3$ and $o\in M$ a
point. Given any conformal class $C\in {\mathcal C}(M)$ and small
$\epsilon>0$ there exists a conformal class $\tilde{C}\in {\mathcal
C}(M)$ such that
$$
\{
\begin{array}{l}
|Y_{\tilde{C}}(M)- Y_C(M)|\leq \epsilon,
\\
|W_{\tilde{C}}(M) -W_C(M)| \leq \epsilon,
\end{array}
\right.
$$
and $\tilde{C}$ is conformally flat near the point $o\in M$.
\end{Proposition}
\begin{Proof}
Let $g\in C$ be any metric. Then following the construction in
\cite[Section 3]{Kobayashi}, we define a smooth family of approximation
metrics.  For a given $o\in M$, there exists a metric $\bar{g}\in
\Riem(M)$, which is conformally flat near the point $o$, such that
$$
\{
\begin{array}{l}
j^1_o(\bar{g})= j^1_o(g),
\\
R_{\bar{g}}(o)=R_g(o).
\end{array}
\right.
$$
For a small $\delta>0$, there exist a positive constant
$0<\epsilon(\delta)$ ($0<\epsilon(\delta)<\delta$) and a
smooth cut-off function $w_{\delta}=w_{\delta}(r)$ ($r\geq 0$) satisfying
$$
0\leq w_{\delta}\leq 1, \ \ \ 
w_{\delta} =\{
\begin{array}{ll}
0 &\mbox{if} \ \ r\geq \delta,
\\
1 &\mbox{if} \ \ 0\leq r\leq \epsilon(\delta),
\end{array}
\right.\ \ \mbox{and} \ \ 
\{
\begin{array}{l}
|r\cdot \dot{w}_\delta| <\delta,
\\
|r^2\cdot \ddot{w}_\delta| <\delta.
\end{array}
\right.
$$
Then let $g_{\delta}:= g + w_{\delta}(\bar{g}-g)$. Then we have that
$$
\{
\begin{array}{ll}
g_{\delta} \to g &\mbox{in the $C^1$-norm on $M$ with respect to $g$},
\\
R_{g_{\delta}} \to R_g & \mbox{in the $C^0$-norm on $M$}
\end{array}
\right.
$$
as $\delta\to 0$. This implies that $Y_{[g_{\delta}]}(M) \to
Y_{[g]}(M)$ as $\delta\to 0$.
\vspace{2mm}

\noindent
To analyze the behaviour of the Weyl constant, we denote by $B_{\delta}$
the disk (with respect to $g$) centered at $o\in M$ of radius
$\delta$. Also we denote $h=\bar{g}-g$ and $T= w_{\delta}\cdot h =
O(r^2)$.  Then
\begin{equation}\label{eq_f1}
(W_{g_{\delta}})^i_{\ jk\ell}- (W_{g})^i_{\ jk\ell} =
O\left(|\nabla^2 T|_g\right) +O\left((1+|T|^2_g)|\nabla T|^2_g
\right).
\end{equation}
Here 
$$
\{
\begin{array}{rcl}
\nabla T &=& \dot{w}_{\delta} \cdot h 
+ w_{\delta}\cdot \nabla h, 
\\
\nabla^2 T &=& \ddot{w}_{\delta}\cdot h + 
2 \dot{w}_{\delta} \cdot \nabla h + w_{\delta}\cdot \nabla^2 h.
\end{array}
\right.
$$
Notice that $\nabla h= O(r)$, $\nabla^2
h= O(1)$, and then
\begin{equation}\label{eq_f2}
\{
\begin{array}{ccl}
|T|^2_g &\leq& K_1 \cdot \delta^4,
\\
|\nabla T |^2_g &\leq& K_2\cdot \delta^2,
\\
|\nabla^2 T |_g &\leq& K_3\cdot \delta + K_4
\end{array}
\right.
\end{equation}
for some positive constants $K_1$, $K_2$, $K_3$, $K_4$. From
(\ref{eq_f1}) and (\ref{eq_f2}), we obtain
$$
|W_{g_{\delta}}|_g \leq |W_g|_g + K_5
$$
for some $K_5>0$ on the disk $B_{\delta}$, and hence there exists a
constant $K_6>0$ such that
$$
|W_{g_{\delta}}|_g^{\frac{n}{2}} \leq
K_6  \left(|W_g|_g^{\frac{n}{2}}+1\right).
$$
This implies that
$$
\begin{array}{rcl}
\displaystyle
\left|
W_{[g_{\delta}]}(M) - W_{[g]}(M)
\right| &=&
\displaystyle
\left|
\int_{B_{\delta}}|W_{g_{\delta}}|_{g_{\delta}}^{\frac{n}{2}} 
d\sigma_{g_{\delta}} -
\int_{B_{\delta}}|W_{g}|_g^{\frac{n}{2}} d\sigma_{g}
\right| 
\\
\\
 &\leq&
\displaystyle
2 K_6 \cdot \int_{B_{\delta}}
\left( 
|W_{g}|_g^{\frac{n}{2}} +1 
\right)d\sigma_{g}
\end{array}
$$
since $d\sigma_{g_{\delta}}\leq 2
d\sigma_{g}$. We obtain that $W_{[g_{\delta}]}(M) \to W_{[g]}(M)$ as
$\delta\to0$.
\end{Proof}
\begin{lproof}{Proof of Theorem A}
First we choose small $\epsilon>0$. Then by definition there exists a
conformal class $C\in {\mathcal C}(M)$ such that
\begin{equation}\label{eq_f3}
\{
\begin{array}{l}
Y_C(M) \geq Y(M) - \frac{\epsilon}{2},
\\
\omega(M)+ \frac{\epsilon}{2} \geq  W_C(M) \geq \omega(M).
\end{array}
\right.
\end{equation}
From Proposition \ref{s3_Prop1} we may assume that $C$ is conformally
flat near some point $o\in M$. 
\vspace{2mm}

\noindent
Now consider the sphere $S^n$.  By an argument similar to the one in
the proof of Theorem \ref{s2_Th1}, for any constant $\hat{\kappa}>0$
there exists a conformal class $\hat{C}\in {\mathcal C}(S^n)$, which
is conformally flat near $\hat{o}\in S^n$, such that
\begin{equation}\label{eq_f4}
\{
\begin{array}{l}
Y_{\hat{C}}(S^n)\geq Y(S^n)- \frac{\epsilon}{2},
\\
\hat{\kappa}+\frac{\epsilon}{2} \geq W_{\hat{C}}(S^n)\geq \hat{\kappa}.
\end{array}
\right.
\end{equation}
We decompose the manifolds $M\setminus \{o\}$ and $S^n\setminus
\{\hat{o}\}$ as follows (see Fig. \ref{s3}.1):
$$
\{
\begin{array}{lcl}
M\setminus \{o\} &= & X\cup (S^{n-1}\times [0,\infty)),
\\
S^n\setminus \{\hat{o}\}&= & \hat{X}\cup (S^{n-1}\times [0,\infty)).
\end{array} \right.  
$$
Let $h_0$ be the standard metric on $S^{n-1}$. We choose metrics
$g\in C$ and $\hat{g}\in \hat{C}$ satisfying
$$
\{
\begin{array}{lcl}
g\in C & \mbox{with} &g(z,t)=h_0(z)+ dt^2 \ \ \mbox{on} \ \ 
S^{n-1}\times [0,\infty),
\\
\hat{g}\in \hat{C} & \mbox{with} &\hat{g}(z,t)=h_0(z)+ dt^2 \ \ \mbox{on} \ \ 
S^{n-1}\times [0,\infty).
\end{array}
\right.
$$
For each $\ell>0$ let $g_{\ell}$ be the natural gluing metric on the
manifold
$$
M\cong M\# S^n  \cong 
\left(X\cup (S^{n-1}\times [0,\ell])\right)\cup_{S^{n-1}\times \{\ell\}}
\left((S^{n-1}\times [0,\ell])\cup \hat{X} \right),
$$
which satisfies 
$$
g_{\ell}|_{X\cup (S^{n-1}\times [0,\ell])} = g ,  \ \ \mbox{and} \ \ 
g_{\ell}|_{\hat{X}\cup (S^{n-1}\times [0,\ell])} = \hat{g}.
$$
\noindent
\hspace*{20mm}\PSbox{sob8.pstex}{5cm}{43mm}
\begin{picture}(0,0)
\put(20,65){{\small $\ell$}}
\put(30,17){{\small $\ell$}}
\put(-140,35){{\small $(M\setminus \{o\}, C)$}}
\put(-110,70){{\small $X$}}
\put(140,35){{\small $\hat{X}$}}
\put(140,100){{\small $(S^n\setminus \{\hat{o}\},\hat{C})$}}
\end{picture}
\vspace{1mm}

\centerline{{\small {\bf Fig. \ref{s3}.1.} The manifolds $M\setminus \{o\}$ and
$S^n\setminus \{\hat{o}\}$.}}
\vspace{2mm}

\noindent
From the argument in the proof of \cite[Theorem 2]{Kobayashi}
(cf. \cite{AB2}), there exists a large constant
$\ell=\ell(\epsilon)>0$ such that
$$
\begin{array}{rcl}
Y_{[g_{\ell}]}(M) &\geq & \displaystyle
Y_{C\sqcup \hat{C}}((M\setminus\{o\})\sqcup (S^{n}\setminus\{\hat{o}\}))
-\frac{\epsilon}{2}
\\
\\
&= & \displaystyle
Y_{C\sqcup \hat{C}}(M\sqcup S^n) -\frac{\epsilon}{2}
= \min\{Y_C(M), Y_{\hat{C}}(S^n)\} -\frac{\epsilon}{2}
\\
\\
&= & \displaystyle
\left(Y(M)-\frac{\epsilon}{2}\right) - \frac{\epsilon}{2} 
= Y(M) -\epsilon.
\end{array}
$$
Recall that $g$ and $\hat{g}$ are conformally flat on
$S^{n-1}\times[0,\infty)$.  Hence for the Weyl constant, we have:
$$
\begin{array}{rcl}
W_{[g_{\ell}]}(M) &=& \displaystyle
\int_{X} |W_g|^{\frac{n}{2}}_g d\sigma_g + 
\int_{\hat{X}} |W_{\hat{g}}|^{\frac{n}{2}}_{\hat{g}} d\sigma_{\hat{g}}
\\
\\
&=& \displaystyle
W_C(M)+ W_{\hat{C}}(S^n).
\end{array}
$$
This combined with (\ref{eq_f3}) and (\ref{eq_f4}) implies
$$
\hat{\kappa} +\omega(M) +\epsilon \geq W_C(M)+ W_{\hat{C}}(S^n)  
\geq \hat{\kappa} +\omega(M). 
$$
We take $\hat{\kappa}=\kappa-\omega(M)>0$ and then obtain the second
inequality in Theorem A. This completes the proof of Theorem A.
\end{lproof}


\vspace{10mm}

\begin{small}
{\sf
\noindent
Kazuo Akutagawa, Shizuoka University, Shizuoka, Japan
\\
e-mail: smkacta\atsign ipc.shizuoka.ac.jp
\\
\\
Boris Botvinnik, University of Oregon, Eugene, USA
\\
e-mail: 
botvinn\atsign math.uoregon.edu
\\
\\
Osamu Kobayashi, Kumamoto University, Kumamoto, Japan
\\
e-mail:
ok\atsign math.sci.kumamoto-u.ac.jp
\\
\\
Harish Seshadri, University of Oregon, Eugene, USA
\\
e-mail: 
harish\atsign darkwing.uoregon.edu 
}
\end{small}
\end{document}